\documentclass[12pt,a4paper,reqno]{amsart}
\usepackage{graphics,epic}
\usepackage{amsmath,amssymb, amsthm}
\usepackage[all,2cell]{xy}

\newtheorem{theorem}{Theorem}[section]
\newtheorem*{theorem*}{Theorem}

\newtheorem{lemma}[theorem]{Lemma}
\newtheorem{proposition}[theorem]{Proposition}
\newtheorem{corollary}[theorem]{Corollary}

\newtheorem*{conjecture*}{Conjecture}

\newtheorem{remark}[theorem]{Remark}

\renewcommand{\hat}[1]{\widehat{#1}}

\newcommand{\ch}{{\rm ch}}
\newcommand{\im}{{\rm im}}
\newcommand{\Hom}{{\rm Hom}\,}
\newcommand{\Glob}{{\rm Glob}\,}

\newcommand{\qdim}{{\rm qdim}\,}

\newcommand{\Z}{\mathbb{Z}}

\newcommand{\C}{\mathbb{C}}

\def\wt{{\rm wt}}
\def\C{{\mathbb C}}

\def\Z{{\mathbb Z}}

\def\1{{\bf 1}}
\def\tr{{\rm tr}}

\def \Hom{{\rm Hom}}

\def \Ind{{\rm Ind}}

\def \pf{\noindent {\bf Proof: \,}}
\def\theequation{5.\arabic{equation}}
\def \h{\mathfrak{h}}
\def \w{\omega}
\def \g{\mathfrak{g}}

\begin{document}

\title[Trace functions and fusion rules of  diagonal coset VOAs]{Trace functions and fusion rules of  diagonal coset vertex operator algebras}

\author{Xingjun Lin}
\address{Xingjun Lin,  School of Mathematics and Statistics, Wuhan University, Wuhan 430072, China.}
\thanks{X. Lin was supported by China NSF grant
11801419 and the starting research fund from Wuhan University (No. 413000076)}
\email{linxingjun88@126.com}
\begin{abstract}
In this paper, irreducible modules of the diagonal coset vertex operator algebra $C(L_{\g}(k+l,0),L_{\g}(k,0)\otimes L_{\g}(l,0))$ are classified  under the assumption that $C(L_{\g}(k+l,0),L_{\g}(k,0)\otimes L_{\g}(l,0))$ is rational, $C_2$-cofinite and certain additional assumption. An explicit  modular transformation formula of traces functions of $C(L_{\g}(k+l,0),L_{\g}(k,0)\otimes L_{\g}(l,0))$ is obtained. As an application, the fusion rules of $C(L_{E_8}(k+2,0), L_{E_8}(k,0)\otimes L_{E_8}(2,0))$ are  determined by using the Verlinde formula.
\end{abstract}
\keywords{Vertex operator algebra; Coset vertex operator algebra; Traces function; Fusion rule; Affine Lie algebra}
\maketitle
\section{Introduction \label{intro}}
\def\theequation{1.\arabic{equation}}
\setcounter{equation}{0}
This paper is a continuation of \cite{Lin} on the diagonal coset vertex operator algebra. Let $\g$ be a finite dimensional simple Lie algebra and $k,l$ be positive integers. Then the irreducible quotient $L_{\g}(k, 0)$ of the vacuum module $V_{\g}(k, 0)$ of the affine Lie algebra $\hat \g$ is a vertex operator algebra \cite{FZ}. By the diagonal action of $\hat \g$, the tensor product vertex operator algebra $L_{\g}(k,0)\otimes L_{\g}(l,0)$ is an integrable $\hat \g$-module  of level $k+l$ \cite{KW}. We use
$C(L_{\g}(k+l,0),L_{\g}(k,0)\otimes L_{\g}(l,0))$ to denote the multiplicity space of $L_{\g}(k+l,0)$ in $L_{\g}(k,0)\otimes L_{\g}(l,0)$.
Then it is known \cite{LL} that $C(L_{\g}(k+l,0),L_{\g}(k,0)\otimes L_{\g}(l,0))$ is a vertex operator algebra, which is called the diagonal coset vertex operator algebra.

One of important problems about the diagonal coset vertex operator algebra is to classify irreducible modules of $C(L_{\g}(k+l,0),L_{\g}(k,0)\otimes L_{\g}(l,0))$. There are several cases that irreducible modules of $C(L_{\g}(k+l,0),L_{\g}(k,0)\otimes L_{\g}(l,0))$ have been classified \cite{JL1}, \cite{JL2}, \cite{ACL}, \cite{CFL}, \cite{Lin}. However, in the general case that $\g$ is a finite dimensional simple Lie algebra and $k,l$ are positive integers, the classification of irreducible modules of $C(L_{\g}(k+l,0),L_{\g}(k,0)\otimes L_{\g}(l,0))$ has not been completed.

 Our first main result is about the classification of irreducible modules of $C(L_{\g}(k+l,0),L_{\g}(k,0)\otimes L_{\g}(l,0))$. Let $L_{\g}(k, \dot{\Lambda})$ (resp. $L_{\g}(l, \ddot{\Lambda})$) be an integrable $\hat \g$-module of level $k$ (reps. $l$).  By the diagonal action of $\hat \g$,  $L_{\g}(k, \dot{\Lambda})\otimes L_{\g}(l, \ddot{\Lambda})$ is an integrable $\hat \g$-module of level $k+l$ \cite{KW}. We use $M_{\dot{\Lambda}, \ddot{\Lambda}}^{\Lambda}$ to denote the multiplicity space of $L_{\g}(k+l, \Lambda)$ in $L_{\g}(k, \dot{\Lambda})\otimes L_{\g}(l, \ddot{\Lambda})$. Then it is known \cite{Lin} that $M_{\dot{\Lambda}, \ddot{\Lambda}}^{\Lambda}$ is a $C(L_{\g}(k+l,0),L_{\g}(k,0)\otimes L_{\g}(l,0))$-module.
 In general, $M_{\dot{\Lambda}, \ddot{\Lambda}}^{\Lambda}$ may be a reducible $C(L_{\g}(k+l,0),L_{\g}(k,0)\otimes L_{\g}(l,0))$-module.
Under the assumption that $C(L_{\g}(k+l,0),L_{\g}(k,0)\otimes L_{\g}(l,0))$ is rational and $C_2$-cofinite, the global dimension of $C(L_{\g}(k+l,0),L_{\g}(k,0)\otimes L_{\g}(l,0))$  and the quantum dimension of $M_{\dot{\Lambda}, \ddot{\Lambda}}^{\Lambda}$ have been obtained in \cite{Lin}. Based on these results, we obtain a sufficient condition such that $M_{\dot{\Lambda}, \ddot{\Lambda}}^{\Lambda}$ is an irreducible $C(L_{\g}(k+l,0),L_{\g}(k,0)\otimes L_{\g}(l,0))$-module (see Theorem \ref{classi}). We further classify irreducible modules of $C(L_{\g}(k+l,0),L_{\g}(k,0)\otimes L_{\g}(l,0))$ under this additional condition (see Theorem \ref{classi}).

 Our second main result is about the modular transformation formula of traces functions of $C(L_{\g}(k+l,0),L_{\g}(k,0)\otimes L_{\g}(l,0))$. We assume that $C(L_{\g}(k+l,0),L_{\g}(k,0)\otimes L_{\g}(l,0))$ is rational and $C_2$-cofinite, then the space spanned by trace functions $Z_M(v, \tau)$ of the irreducible $C(L_{\g}(k+l,0),L_{\g}(k,0)\otimes L_{\g}(l,0))$-modules  $M$ affords a representation of the modular group $SL(2, \Z)$ \cite{Z}. Our goal is to determine this representation explicitly. In the case that $v$ is the vacuum vector of $C(L_{\g}(k+l,0),L_{\g}(k,0)\otimes L_{\g}(l,0))$,  trace functions $Z_M(v, \tau)$ are a special kind of branching functions studied in \cite{KW}. An explicit modular transformation formula for branching functions has been obtained in \cite{KW}. In this paper, we show that the same transformation formula for branching functions is valid for the trace functions (see Theorem \ref{main}).  The main idea of the proof is similar to that in \cite{DKR}, in which the modular transformation formula of trace functions of  parafermion vertex operator algebras has been obtained.

As an application,  the fusion rules of  $C(L_{E_8}(k+2,0),L_{E_8}(k,0)\otimes L_{E_8}(2,0))$ are determined. It has been proved in \cite{Lin} that $C(L_{E_8}(k+2,0),L_{E_8}(k,0)\otimes L_{E_8}(2,0))$ is rational and $C_2$-cofinite. Moreover, irreducible modules of $C(L_{E_8}(k+2,0),L_{E_8}(k,0)\otimes L_{E_8}(2,0))$ have been classified in \cite{Lin}. Therefore, we can use Theorem \ref{main} to obtain an explicit modular transformation formula of trace functions of  $C(L_{E_8}(k+2,0),L_{E_8}(k,0)\otimes L_{E_8}(2,0))$. As a consequence, the fusion rules of $C(L_{E_8}(k+2,0),L_{E_8}(k,0)\otimes L_{E_8}(2,0))$ are determined in Theorem \ref{fusion}  by using the Verlinde formula.

The paper is organized as follows: In Section 2, we recall some facts about  generalized theta functions of vertex operator algebras, which play an important role in the proof of the modular transformation formula of traces functions of $C(L_{\g}(k+l,0),L_{\g}(k,0)\otimes L_{\g}(l,0))$. In Section 3, we recall some facts about  $C(L_{\g}(k+l,0),L_{\g}(k,0)\otimes L_{\g}(l,0))$-modules. In Section 4, we classify irreducible  $C(L_{\g}(k+l,0),L_{\g}(k,0)\otimes L_{\g}(l,0))$-modules under the assumption that $C(L_{\g}(k+l,0),L_{\g}(k,0)\otimes L_{\g}(l,0))$ is rational, $C_2$-cofinite and certain additional assumption. In Section 5, an explicit modular transformation formula of trace functions of $C(L_{\g}(k+l,0),L_{\g}(k,0)\otimes L_{\g}(l,0))$ is obtained. In Section 6, we determine the fusion rules of  $C(L_{E_8}(k+2,0),L_{E_8}(k,0)\otimes L_{E_8}(2,0))$.

\section{Preliminaries }
\def\theequation{2.\arabic{equation}}
\setcounter{equation}{0}
In this section, we recall some facts about generalized theta functions of vertex operator algebras, which play an important role in Section 5. We will continue to use notations in \cite{Lin}.
\subsection{Modular invariance properties of trace functions of vertex operator algebras}
In this subsection, we recall from \cite{Z} some facts about  modular invariance properties of trace functions of vertex operator algebras.  Let $(V, Y(\cdot, z), \1, \w)$ be a rational and $C_2$-cofinite  vertex operator algebra. It was proved in \cite{DLM1} that $V$ has only finitely many inequivalent irreducible admissible modules $M^0, M^1,...,M^p$. Moreover, $M^i, 0\leq i \leq p$, has the form $$M^i=\bigoplus_{n=0}^{\infty}M^i_{\lambda_i+n},$$ with $M^i_{\lambda_i}\neq 0$ for some number $\lambda_i$ which is called {\em conformal weight} of $M^i$. For any irreducible $V$-module $M^i$, the {\em trace function} associated to $M^i$ is defined as follows: For any homogenous element $v\in V$ and $\tau\in \mathcal H=\{\tau\in \mathbb{C}| \im\tau>0\}$,
\begin{equation*}
Z_{M^i}(v,\tau):=\tr_{M^i}o(v)q^{L(0)-c/24}=q^{\lambda_i-c/24}\sum_{n\in\mathbb{Z}_{\geq 0}} \tr_{M^i_{\lambda_i+n}}o(v)q^n,
\end{equation*}
where $o(v)=v_{\wt v-1}$ and $q=e^{2\pi \sqrt{-1}\tau}$. Since $V$ is $C_2$-cofinite, $Z_{M^i}(v,\tau)$ converges to a holomorphic function on the domain $|q| < 1$ \cite{DLM2}, \cite{Z}. The $Z_M(\1, \tau)$ which is also denoted by $\ch_q M$ is called the $q$-character of $M$.

To study modular invariance properties of trace functions of vertex operator algebras, another vertex operator algebra structure $(V, Y[\cdot, z], \1, \w-c/24)$ is defined on $V$ in \cite{Z} such that for any homogenous element $v\in V$,  $$Y[v,z]=\sum_{n\in \Z}v_{[n]}z^{-n-1}=Y(v, e^{z}-1)e^{z\wt v}.$$  Then we have
\begin{align}\label{trans}
v_{[n]}=n!\sum_{i\geq n}c(\wt v, i, n)v_i,
\end{align}
where the numbers $c(\wt v, i, n)$ are determined by $$\left(\begin{tabular}{c}
$\wt v-1+z$\\
$i$ \\
\end{tabular}\right)=\sum_{n=0}^ic(\wt v, i, n)z^n.$$  Write $Y[\w-c/24, z]=\sum_{n\in \Z}L[n]z^{-n-2}$. Then $V$ has a new grading
$$V=\oplus_{n\geq 0} V_{[n]},$$
where $ V_{[n]}=\{v\in V|L[0]v=nv\}$.
For $v\in V_{[n]}$, we set $\wt[v]=n$. Then the following result has been obtained in \cite{Z} (also see \cite{DLM2}).
\begin{theorem}\label{minvariance}
 Let $V$ be a rational and $C_2$-cofinite vertex operator algebra, $M^0,...,M^p$  be all the inequivalent irreducible $V$-modules.   Then there is a group homomorphism $\rho_V: SL(2, \Z)\to GL_{p+1}(\C)$ with $\rho_V(\gamma)=(\gamma_{i,j})$ such that for any $0\leq i\leq p$ and homogeneous $v\in V_{[n]}$, $$Z_{M^i}(v,\frac{a\tau+b}{c\tau+d})=(c\tau+d)^n\sum_{j=0}^p \gamma_{i,j}Z_{M^j}(v,\tau).$$
 Moreover, the matrix $(\gamma_{i,j})$ is independent of the choice of $v\in V$.
\end{theorem}
Recall that the full modular group $SL(2, \mathbb{Z})$ has generators $S=\left(\begin{array}{cc}0 & -1\\ 1 & 0\end{array}\right)$, $T=\left(\begin{array}{cc}1 & 1\\ 0 & 1\end{array}\right)$. It is known that $\rho_V(S)$ is a symmetric matrix. We will also use $S=(S_{i,j})$ to denote the  matrix $\rho_V(S)$.
\subsection{Modular invariance properties of the generalized theta functions} In this subsection, we will recall some facts about the generalized theta functions from \cite{DKR} and \cite{Kra}. Let $V$ be a strongly regular vertex operator algebra, $M^0,...,M^p$  be all the inequivalent irreducible $V$-modules.  Then it was proved in \cite{DM1} that $V_1$ is a reductive Lie algebra. Let $\h$ be a Cartan subalgebra of $V_1$. Then the abelian Lie algebra $\h$ acts on $M^i$ semisimplely for all $i$.

Recall from \cite{FHL} that a bilinear form $(\cdot, \cdot)$ on $V$ is called invariant if
$$( Y(a, z)u, v)= ( u, Y(e^{zL(1)}(-z^{-2})^{L(0)}a, z^{-1})v),$$for $a, u, v\in V$. Since $V$ is a strongly regular vertex operator algebra, it was proved in \cite{Li} that there is a unique nondegenerate invariant bilinear form $(\cdot, \cdot)$ on $V$ such that $( u, v) \1=u_1v$ for $u, v\in V_1$.

Define $$\chi_i(v, q)=\tr_{M^i}e^{2\pi \sqrt{-1} v_0}q^{L(0)-c/24}$$ for $v\in \h$. Then the following results have been proved in \cite{DLMa} (see also \cite{DKR}).
\begin{theorem}
Let $V$ be a strongly regular vertex operator algebra, $M^0,...,M^p$  be all the inequivalent irreducible $V$-modules. Then\\ 
(i) For $v\in \h$, $\chi_i(v, q)$ converges to a homomorphic function $\chi_i(v, \tau)$ in $\mathcal{H}$ with $q=e^{2\pi\sqrt{-1}\tau}$.\\
(ii) For any $\gamma=\left(\begin{array}{cc}a & b\\ c& d\end{array}\right)\in SL(2, \Z)$ and $v\in \h$,
$$\chi_i(\frac{v}{c\tau+d}, \frac{a\tau+d}{c\tau+d})=e^{\pi\sqrt{-1}(c( v, v)/(c\tau+d))}\sum_{j=0}^p\gamma_{i, j}\chi_j(v, \tau),$$ where $\gamma_{i, j}$ is the same as in Theorem \ref{minvariance}.
\end{theorem}

We now let $U$ and $W$ be strongly regular vertex operator algebras such that $U\otimes W$ is a vertex subalgebra of $V$ and that $U\otimes W$, $V$ have the same Virasoro vector. Then $U_1$ is a reductive Lie algebra and we fix a Cartan subalgebra $\h^{U}$ of $U_1$. Following \cite{Kra}, we define
$$\Phi_{j}(w, u,v, q)=\tr_{M^j}o(w)e^{2\pi \sqrt{-1}(o(u)+\frac{( u, v)}{2})}q^{L(0)-c/24+o(v)+\frac{( v, v)}{2}},$$
for $w\in W$, $u, v\in \h^U$. By the formula (\ref{trans}), the following result follows from Theorem 1.1 of \cite{Kra}.
\begin{theorem}\label{theta1}
Let $V$, $U$ and $W$ be strongly regular vertex operator algebras such that $U\otimes W$ is a vertex subalgebra of $V$ and that $U\otimes W$, $V$ have the same Virasoro vector, $M^0,...,M^p$  be all the inequivalent irreducible $V$-modules. Assume that $\Phi_{j}(w, u,v, q)$ converges to a holomorphic function $\Phi_{j}(w, u,v, \tau)$ in $\mathcal{H}$ with $q=e^{2\pi \sqrt{-1}\tau}$ for any homogeneous $w\in W$ and $u,v\in \h^U$. Then for any $\gamma=\left(\begin{array}{cc}a & b\\ c& d\end{array}\right)\in SL(2, \Z)$,
$$\Phi_{i}(w, u,v, \frac{a\tau+b}{c\tau+d})=(c\tau+d)^{\wt [w]}\sum_{j=0}^p \gamma_{i, j}\Phi_{j}(w, bv+du, av+cu,\tau)$$ where $\gamma_{i, j}$ is the same as in Theorem \ref{minvariance}.
\end{theorem}

Set $\chi_{j}(w, u, q)=\Phi_j(w, u, 0, q)$. Then we have $$\chi_{j}(w, u, q)=\tr_{M^j}o(w)e^{2\pi \sqrt{-1}o(u)}q^{L(0)-c/24},$$
for $w\in W$, $u\in \h^U$.  Using Theorem \ref{theta1}, one can show that the following result.
\begin{theorem}\label{theta2}
Let $V$, $U$ and $W$ be strongly regular vertex operator algebras such that $U\otimes W$ is a vertex subalgebra of $V$ and that $U\otimes W$, $V$ have the same Virasoro vector, $M^0,...,M^p$  be all the inequivalent irreducible $V$-modules. Assume that $\chi_{j}(w, u, q)$ converges to a holomorphic function $\chi_{j}(w, u, \tau)$ in $\mathcal{H}$ with $q=e^{2\pi \sqrt{-1}\tau}$ for any homogeneous $w\in W$ and $u\in \h^U$. Then for any $\gamma=\left(\begin{array}{cc}a & b\\ c& d\end{array}\right)\in SL(2, \Z)$,
$$\chi_{i}(w, \frac{u}{c\tau+d}, \frac{a\tau+b}{c\tau+d})=(c\tau+d)^{\wt[w]} e^{\pi\sqrt{-1}(c( u, u)/(c\tau+d))}\sum_{j=0}^p \gamma_{i, j}\chi_{j}(w, u,\tau)$$ where $\gamma_{i, j}$ is the same as in Theorem \ref{minvariance}.
\end{theorem}

\section{Modules of  diagonal coset vertex operator algebras}
\def\theequation{3.\arabic{equation}}
\setcounter{equation}{0}
\subsection{Diagonal coset vertex operator algebras}In this subsection, we shall recall some facts about diagonal coset vertex operator algebras from \cite{Lin}. Let $\g$ be a finite dimensional simple Lie algebra and $\langle\, ,\, \rangle$ the normalized Killing form of $\g$, i.e., $\langle\theta, \theta\rangle=2$ for the highest root $\theta$ of $\g$. The affine Lie algebra associated to $\g$ is defined on $\hat{\g}=\g\otimes \C[t^{-1}, t]\oplus \C K$ with Lie brackets
\begin{align*}
[x(m), y(n)]&=[x, y](m+n)+\langle x, y\rangle m\delta_{m+n,0}K,\\
[K, \hat\g]&=0,
\end{align*}
for $x, y\in \g$ and $m,n \in \Z$, where $x(n)$ denotes $x\otimes t^n$.

For a positive integer $k$, define
\begin{align*}
V_{\g}(k, 0)=\Ind_{\g\otimes \C[t]\oplus \C K}^{\hat \g}\C,
\end{align*}
where $\C$ is viewed as a module of $\g\otimes \C[t]\oplus \C K$ such that $\g\otimes \C[t]$ acts as $0$ and $K$ acts as $k$. It is well-known that $V_{\g}(k, 0)$ has a unique maximal proper submodule which is denoted by $J(k, 0)$ (see \cite{K}). Let $L_{\g}(k, 0)$ be the corresponding irreducible quotient module. It was proved in \cite{FZ} that $L_{\g}(k, 0)$ has a vertex operator algebra structure such that the Virasoro vector
\begin{align*}
\w=\frac{1}{2(k+h^{\vee})} \sum_{i=1}^{\dim \g} u_i(-1)u_i(-1)\1,
\end{align*}
where $h^\vee$ denotes the dual Coxeter number of $\g$ and $\{u_i|1\leq i\leq \dim \g\}$ is an orthonormal basis of $\g$ with respect to $\langle,\rangle$. Moreover, the following result has been proved in \cite{DLM2}, \cite{FZ}.
\begin{theorem}
Let $k$ be a positive integer. Then $L_{\g}(k, 0)$ is a strongly regular vertex operator algebra.
\end{theorem}

 For positive integers  $k$, $l$, we consider the vertex operator algebra $L_{\g}(k,0)\otimes L_{\g}(l,0)$. Let $W$ be the vertex subalgebra of  $L_{\g}(k,0)\otimes L_{\g}(l,0)$ generated by $\{x(-1)\1\otimes \1+\1\otimes x(-1)\1|x\in \g\}$. By Theorem 3.1 of \cite{DM2}, $W$ is isomorphic to $L_{\g}(k+l,0)$. Set \begin{align*}
C(L_{\g}&(k+l,0),L_{\g}(k,0)\otimes L_{\g}(l,0))
\\&=\{u\in L_{\g}(k,0)\otimes L_{\g}(l,0)|u_{n}v=0, \forall v\in L_{\g}(k+l,0), \forall n\in \Z_{\geq 0}\}.
\end{align*}
It is well-known that $C(L_{\g}(k+l,0),L_{\g}(k,0)\otimes L_{\g}(l,0))$ is a vertex subalgebra of $L_{\g}(k,0)\otimes L_{\g}(l,0)$ \cite{LL}. Moreover, we have the following result (see Lemma 3.5 of \cite{Lin}).
\begin{proposition}\label{scoset}
Let $\w^1, \w^2, \w^a$ be the Virasoro vectors of $L_{\g}(k,0)$, $L_{\g}(l,0)$, $L_{\g}(k+l,0)$, respectively. Then $C(L_{\g}(k+l,0),L_{\g}(k,0)\otimes L_{\g}(l,0))$  is a  simple vertex operator algebra with the Virasoro vector $\w^1+\w^2-\w^a$.
\end{proposition}

\subsection{Modules of diagonal coset vertex operator algebras}
 In this subsection, we recall from \cite{Lin} some facts about modules of the vertex operator algebra $C(L_{\g}(k+l,0),L_{\g}(k,0)\otimes L_{\g}(l,0))$. Let $\g$ be a finite dimensional simple Lie algebra and $\h$ be a Cartan subalgebra of $\g$. We denote the corresponding root system by $\Delta_{\g}$ and the root lattice by $Q$. Then the weight lattice $P$ of $\g$ is the set of $\lambda\in \h$ such that $\frac{2\langle\lambda, \alpha\rangle}{\langle\alpha, \alpha\rangle}\in\Z$ for all $\alpha\in \Delta_{\g}$. We choose a set $\{\alpha_1,\cdots,\alpha_l\}$  of simple roots, and denote the set of positive roots by $\Delta_{\g}^+$.  Note that $P$ is equal to $\oplus_{i=1}^l\Z\Lambda_i$, where $\Lambda_i$ are the fundamental weights defined by $\frac{2\langle\Lambda_i, \alpha_j\rangle}{\langle\alpha_j, \alpha_j\rangle}=\delta_{i,j}$. We also use the standard notation $P_+$ to denote the set $\{\Lambda\in P\mid\frac{2\langle\Lambda, \alpha_j\rangle}{\langle\alpha_j, \alpha_j\rangle}\geq 0,~1\leq j\leq l \}$ of dominant weights.

For a positive integer $k$ and a weight $\Lambda \in P$, we use  $L_{\g}(\Lambda)$ to denote the irreducible highest weight module of $\g$ with highest weight $\Lambda$. We then define
\begin{align*}
V_{\g}(k, \Lambda)=\Ind_{\g\otimes \C[t]\oplus \C K}^{\hat \g}L_{\g}(\Lambda),
\end{align*}
where $L_{\g}(\Lambda)$ is viewed as a module of $\g\otimes \C[t]\oplus \C K$ such that $\g\otimes t\C[t]$ acts as $0$ and $K$ acts as $k$. It is well-known that $V_{\g}(k, \Lambda)$ has a unique maximal proper submodule which is denoted by $J(k, \Lambda)$ (see \cite{K}). Let $L_{\g}(k, \Lambda)$ be the corresponding irreducible quotient module. Then the following results have been proved in  \cite{FZ}, \cite{K}.
\begin{theorem}\label{moduleaff}
Let $k$ be a positive integer. Then\\
 (1) $L_{\g}(k, \Lambda)$ is a module for the vertex operator algebra $L_{\g}(k, 0)$ if and only if $\Lambda \in P_+^k$, where $P_+^k=\{\Lambda \in P_+|\langle\Lambda, \theta\rangle\leq k\}$.\\
 (2) If $L_{\g}(k, \Lambda)$ is an $L_{\g}(k, 0)$-module  such that $L_{\g}(k, \Lambda)\ncong L_{\g}(k, 0)$, then the conformal weight of $L_{\g}(k, \Lambda)$ is positive.
\end{theorem}

 We now recall some facts about modules of the vertex operator algebra $C(L_{\g}(k+l,0),L_{\g}(k,0)\otimes L_{\g}(l,0))$. For $\dot\Lambda\in P_+^{k}, \ddot\Lambda\in P_+^{l}$, it follows from Theorem \ref{moduleaff} that $L_{\g}(k, \dot{\Lambda})$ and $L_{\g}(l, \ddot{\Lambda})$ are $L_{\g}(k, 0)$-module and $L_{\g}(l, 0)$-module, respectively. Then $L_{\g}(k, \dot{\Lambda})\otimes L_{\g}(l, \ddot{\Lambda})$ is an $L_{\g}(k, 0)\otimes L_{\g}(l, 0)$-module. As a consequence, $L_{\g}(k, \dot{\Lambda})\otimes L_{\g}(l, \ddot{\Lambda})$ may be viewed as an $L_{\g}(k+l, 0)$-module. Since $L_{\g}(k+l, 0)$ is strongly regular, $L_{\g}(k, \dot{\Lambda})\otimes L_{\g}(l, \ddot{\Lambda})$ is completely reducible as an $L_{\g}(k+l, 0)$-module. For any $\Lambda\in P_+^{k+l}$, we define $$M_{\dot{\Lambda}, \ddot{\Lambda}}^{\Lambda}=\Hom_{L_{\g}(k+l, 0)} (L_{\g}(k+l, \Lambda), L_{\g}(k, \dot{\Lambda})\otimes L_{\g}(l, \ddot{\Lambda})).$$ Then it was proved in \cite{Lin} that $C(L_{\g}(k+l,0),L_{\g}(k,0)\otimes L_{\g}(l,0))=M_{0, 0}^0$ and $M_{\dot{\Lambda}, \ddot{\Lambda}}^{\Lambda}$ is a $C(L_{\g}(k+l,0),L_{\g}(k,0)\otimes L_{\g}(l,0))$-module. Moreover, we have
  the following results, which were essentially established in \cite{KW} (see also \cite{Lin}).
\begin{proposition}\label{decomp}
Let $\theta=\sum_{i=1}^l a_i\alpha_i$, $a_i\in \Z_+$, be the highest root of $\g$ and set $J=\{i|a_i=1\}$. Then for $\dot\Lambda\in P_+^{k}, \ddot\Lambda\in P_+^{l}, \Lambda\in P_+^{k+l}$, we have\\
(1) $M_{\dot{\Lambda}, \ddot{\Lambda}}^{\Lambda}\neq 0$ if and only if $\dot{\Lambda}+\ddot{\Lambda}-\Lambda\in Q$.\\
(2) $L_{\g}(k, \dot{\Lambda})\otimes L_{\g}(l, \ddot{\Lambda})$ viewed as an $L_{\g}(k+l, 0)\otimes C(L_{\g}(k+l,0),L_{\g}(k,0)\otimes L_{\g}(l,0))$-module has the following decomposition
$$L_{\g}(k, \dot{\Lambda})\otimes L_{\g}(l, \ddot{\Lambda})=\oplus_{\Lambda \in P_+^{k+l}; \dot{\Lambda}+\ddot{\Lambda}-\Lambda\in Q}L_{\g}(k+l, \Lambda)\otimes M_{\dot{\Lambda}, \ddot{\Lambda}}^{\Lambda}.$$
(3) Viewed as a $C(L_{\g}(k+l,0),L_{\g}(k,0)\otimes L_{\g}(l,0))$-module, the conformal weight of $M_{\dot{\Lambda}, \ddot{\Lambda}}^{\Lambda}$ is equal to or larger than $0$. Moreover, the conformal weight of $M_{\dot{\Lambda}, \ddot{\Lambda}}^{\Lambda}$ is equal to $0$ only if $(\dot{\Lambda}, \ddot{\Lambda},\Lambda)=(k\Lambda_i, l\Lambda_i,(k+l)\Lambda_i), i\in J$, where $\Lambda_0=0$.
\end{proposition}

We next show that there may be isomorphisms between $C(L_{\g}(k+l,0),L_{\g}(k,0)\otimes L_{\g}(l,0))$-modules $\{M_{\dot{\Lambda}, \ddot{\Lambda}}^{\Lambda}|\dot\Lambda\in P_+^{k}, \ddot\Lambda\in P_+^{l}, \Lambda\in P_+^{k+l}\}$.  For any $h\in \h$, set
$$\Delta(h, z)=z^{h(0)}\exp\left(\sum_{n=1}^{\infty}\frac{h(n)(-z)^{-n}}{-n}\right).$$
Then the following result has been established in Proposition 5.4 of \cite{L4}.
\begin{proposition}
Let $M$ be an irreducible $L_{\g}(k, 0)$-module and $h$ be an element of $\h$ such that $h(0)$ has only integral eigenvalues on $L_{\g}(k, 0)$. Set
$$(M^{(h)}, Y_{M^{(h)}}(\cdot, z))=(M, Y(\Delta(h, z)\cdot, z)).$$ Then $(M^{(h)}, Y_{M^{(h)}}(\cdot, z))$ is an irreducible  $L_{\g}(k, 0)$-module.
\end{proposition}
 Let $h^i\in \h$ for $i=1, \cdots, l$ defined by $\alpha_i(h^j)=\delta_{i,j}$ for $j=1, \cdots, l$. Set $$P^{\vee}=\Z h^1+\cdots+\Z h^l.$$ Then it is known \cite{L3} that $h(0)$ has only integral eigenvalues on $L_{\g}(k, 0)$ if and only if $h\in P^{\vee}$. In particular, for any irreducible $L_{\g}(k, 0)$-module $M$ and $h\in P^{\vee}$, $M^{(h)}$ is also an irreducible $L_{\g}(k, 0)$-module.  We now let $\{\alpha_1^{\vee},\cdots,\alpha_l^{\vee}\}$ be the set of simple coroots and $Q^{\vee}=\Z\alpha_1^{\vee}+\cdots+\Z\alpha_l^{\vee}$. Then it was proved in Proposition 2.25 of \cite{L3} that $M^{(h)}\cong M$ for any irreducible $L_{\g}(k, 0)$-module $M$ and $h\in Q^{\vee}$. Therefore, this induces  an action of $P^{\vee}/Q^{\vee}$ on the set $\{L_{\g}(k, \Lambda)|\Lambda \in P_+^k\}$ (see Proposition 2.24 of \cite{L3}). Moreover, the following result has been established in Theorem 2.26 of \cite{L3}.
 \begin{proposition}\label{iden1}
 (1) $|P^{\vee}/Q^{\vee}|=|P/Q|=|J|+1$.\\
 (2) $P^{\vee}/Q^{\vee}=\{ 0+Q^{\vee}\}\cup\{h^i+Q^{\vee}|i\in J\}$.
 \end{proposition}

Thus, for any $h+Q^{\vee}\in P^{\vee}/Q^{\vee}$, there exists an element $\Lambda^{(h)}\in P_+^k$ such that  $L_{\g}(k, \Lambda)^{(h)}$ is isomorphic to $L_{\g}(k, \Lambda^{(h)})$.  As a result, we obtain the following isomorphisms between $C(L_{\g}(k+l,0),L_{\g}(k,0)\otimes L_{\g}(l,0))$-modules by using the operator $\Delta(\cdot, z)$ (see Corollary 4.4 of \cite{Lin}).
\begin{proposition}\label{iden}
For any  $h+Q^{\vee}\in P^{\vee}/Q^{\vee}$ and $\dot\Lambda\in P_+^{k}, \ddot\Lambda\in P_+^{l}, \Lambda\in P_+^{k+l}$, we have $M_{\dot{\Lambda}, \ddot{\Lambda}}^{\Lambda}\cong M_{\dot{\Lambda}^{(h)}, \ddot{\Lambda}^{(h)}}^{\Lambda^{(h)}}$ as $C(L_{\g}(k+l,0),L_{\g}(k,0)\otimes L_{\g}(l,0))$-modules.
\end{proposition}

\section{Classification of irreducible modules of  diagonal coset vertex operator algebras}
\def\theequation{4.\arabic{equation}}
\setcounter{equation}{0}
In this section, we shall classify irreducible modules of $C(L_{\g}(k+l,0),L_{\g}(k,0)\otimes L_{\g}(l,0))$ under the assumption that $C(L_{\g}(k+l,0),L_{\g}(k,0)\otimes L_{\g}(l,0))$ is rational, $C_2$-cofinite and certain additional assumption.
\subsection{Quantum dimensions of  modules of  diagonal coset vertex operator algebras}
In this subsection,  we recall some facts about quantum dimensions of  modules of  diagonal coset vertex operator algebras from \cite{Lin}. Let $V$ be a strongly regular vertex operator algebra and $M^0=V, M^1,...,M^p$ be all the inequivalent irreducible $V$-modules.  For a $V$-module $M$, the {\em quantum dimension} of $M$ is defined to be
\begin{align*}
\qdim_V M=\lim_{y\to 0^+}\frac{Z_{M}(\1, \sqrt{-1}y)}{Z_{V}(\1, \sqrt{-1}y)},
\end{align*}
where $y$ is real and positive. The {\em global dimension}  of $V$ is defined to be
\begin{align*}
{\rm Glob}\, V=\sum_{i=0}^p(\qdim_V M^i)^2.
\end{align*}
The following results were proved in  \cite{DJX}.
\begin{theorem}\label{qdim1}
Let $V$ be a strongly regular vertex operator algebra and $M^0=V, M^1,...,M^p$ be all the inequivalent irreducible $V$-modules. Assume that the conformal weights of $ M^1,...,M^p$ are greater than $0$. Then \\
(1) $\qdim M^i\geq 1$ for any $0\leq i\leq p$.\\
(2) $\qdim_V M^i=\frac{S_{0,i}}{S_{0,0}}$.\\
(3) ${\rm Glob}\, V=\frac{1}{S_{0,0}^2}.$
\end{theorem}

To determine quantum dimensions of  modules of  diagonal coset vertex operator algebras, we need to recall some facts about modular invariance properties of affine vertex operator algebras. Let $\g$ be a finite dimensional simple Lie algebra and $k$ be a positive integer. Then $L_{\g}(k, 0)$ is strongly regular and $\{L_{\g}(k, \Lambda)|\Lambda \in P_+^k\}$ are all the inequivalent irreducible $L_{\g}(k, 0)$-modules. Hence, the set $\{Z_{L_{\g}(k, \Lambda)}(v,\tau)|\Lambda \in P_+^k\}$  of trace functions is closed under the action of $SL(2,\Z)$. In particular, the following result was essentially obtained in \cite{K,KW}.
\begin{theorem}\label{miden}
 For $\Lambda, \Lambda' \in P_+^k$, let $S_{L_{\g}(k, \Lambda), L_{\g}(k, \Lambda')}$ be complex numbers such that
 $$Z_{L_{\g}(k, \Lambda)}(v,\frac{-1}{\tau})=\tau^{\wt [v]}\sum_{\Lambda' \in P_+^k}S_{L_{\g}(k, \Lambda), L_{\g}(k, \Lambda')}Z_{L_{\g}(k, \Lambda')}(v,\tau)$$
 for any homogeneous vector $v\in L_{\g}(k, 0)$.
 Then
$$S_{L_{\g}(k, 0), L_{\g}(k, \Lambda)}=|P/(k+h^\vee)Q_L|^{-1/2} (k+h^\vee)^{-l/2}\prod_{\alpha\in \Delta_{\g}^+}2\sin \frac{\pi \langle\Lambda +\rho, \alpha\rangle}{k+h^\vee},$$
where $Q_L\subseteqq Q$ is the sublattice of $Q$ spanned by all long roots and $\rho=\sum_{i=1}^l\Lambda_i$. Moreover, $\sum_{\Lambda\in P^k_+} S_{L_{\g}(k, 0), L_{\g}(k, \Lambda)}^2=1$.
\end{theorem}

To determine the global dimension of $C(L_{\g}(k+l,0),L_{\g}(k,0)\otimes L_{\g}(l,0))$, we assume that the vertex operator algebra $C(L_{\g}(k+l,0),L_{\g}(k,0)\otimes L_{\g}(l,0))$ is rational and $C_2$-cofinite. Then the following results have been proved in Theorem 4.6 of \cite{Lin}.
\begin{theorem}\label{crational}
Let $k, l$ be positive integers. Suppose that the vertex operator algebra $C(L_{\g}(k+l,0),L_{\g}(k,0)\otimes L_{\g}(l,0))$ is rational and $C_2$-cofinite. Then\\
(1) $C(L_{\g}(k+l,0),L_{\g}(k,0)\otimes L_{\g}(l,0))$ is strongly regular. \\
(2) Any irreducible $C(L_{\g}(k+l,0),L_{\g}(k,0)\otimes L_{\g}(l,0))$-module is isomorphic to a submodule of $M_{\dot{\Lambda}, \ddot{\Lambda}}^{\Lambda}$ for some $\dot\Lambda\in P_+^{k}, \ddot\Lambda\in P_+^{l}, \Lambda\in P_+^{k+l}$.\\
(3) All the conformal weights of irreducible $C(L_{\g}(k+l,0),L_{\g}(k,0)\otimes L_{\g}(l,0))$-modules except $C(L_{\g}(k+l,0),L_{\g}(k,0)\otimes L_{\g}(l,0))$ are larger than $0$.
\end{theorem}

Under the assumption that $C(L_{\g}(k+l,0),L_{\g}(k,0)\otimes L_{\g}(l,0))$ is rational and $C_2$-cofinite, the global dimension of $C(L_{\g}(k+l,0),L_{\g}(k,0)\otimes L_{\g}(l,0))$ has been determined in Theorem 4.7 of \cite{Lin}.
\begin{theorem}\label{qdimc}
Let $k, l$ be positive integers. Suppose that the vertex operator algebra $C(L_{\g}(k+l,0),L_{\g}(k,0)\otimes L_{\g}(l,0))$ is rational and $C_2$-cofinite. Then we have\\
(1) $\Glob C(L_{\g}(k+l,0),L_{\g}(k,0)\otimes L_{\g}(l,0))=\frac{1}{|P/Q|^2S_{L_{\g}(k, 0), L_{\g}(k, 0)}^2S_{L_{\g}(l, 0), L_{\g}(l, 0)}^2S_{L_{\g}(k+l, 0), L_{\g}(k+l, 0)}^2 }$.\\
(2) For any $\dot\Lambda\in P_+^{k}, \ddot\Lambda\in P_+^{l}, \Lambda\in P_+^{k+l}$ such that $\dot{\Lambda}+\ddot{\Lambda}-\Lambda\in Q$, we have $$\qdim_{C(L_{\g}(k+l,0),L_{\g}(k,0)\otimes L_{\g}(l,0))} M_{\dot{\Lambda}, \ddot{\Lambda}}^{\Lambda}=\frac{S_{L_{\g}(k, 0), L_{\g}(k, \dot\Lambda)}S_{L_{\g}(l, 0), L_{\g}(l, \ddot\Lambda)}S_{L_{\g}(k+l, 0), L_{\g}(k+l, \Lambda)}}{S_{L_{\g}(k, 0), L_{\g}(k, 0)}S_{L_{\g}(l, 0), L_{\g}(l, 0)}S_{L_{\g}(k+l, 0), L_{\g}(k+l, 0)} }.$$
\end{theorem}
\subsection{Classification of irreducible modules of $C(L_{\g}(k+l,0),L_{\g}(k,0)\otimes L_{\g}(l,0))$} Let $\g$ be a finite dimensional simple Lie algebra and  $k, l$ be positive integers. In this subsection, we shall classify irreducible modules of $C(L_{\g}(k+l,0),L_{\g}(k,0)\otimes L_{\g}(l,0))$ under certain assumption. First, we have the following results, which have been proved in Theorem 5.1 of \cite{Lin}.
\begin{theorem}\label{general}
Let $\g$ be a finite dimensional simple Lie algebra and  $k, l$ be positive integers. Suppose the following two conditions hold:\\
(i) The vertex operator algebra $C(L_{\g}(k+l,0),L_{\g}(k,0)\otimes L_{\g}(l,0))$ is rational and $C_2$-cofinite.\\
(ii)  There exist $C(L_{\g}(k+l,0),L_{\g}(k,0)\otimes L_{\g}(l,0))$-modules $W^1, \cdots, W^s$ such that for any $\dot\Lambda\in P_+^{k}, \ddot\Lambda\in P_+^{l}, \Lambda\in P_+^{k+l}$, $ M_{\dot{\Lambda}, \ddot{\Lambda}}^{\Lambda}$ is a direct sum of some $W^{i_1}, \cdots, W^{i_t}$. Moreover, $$\sum_{i=1}^s (\qdim_{C(L_{\g}(k+l,0),L_{\g}(k,0)\otimes L_{\g}(l,0))} W^i)^2=\Glob C(L_{\g}(k+l,0),L_{\g}(k,0)\otimes L_{\g}(l,0)).$$
Then $W^1, \cdots, W^s$ are irreducible $C(L_{\g}(k+l,0),L_{\g}(k,0)\otimes L_{\g}(l,0))$-modules. Moreover, $W^1, \cdots, W^s$ are all the inequivalent irreducible $C(L_{\g}(k+l,0),L_{\g}(k,0)\otimes L_{\g}(l,0))$-modules.
\end{theorem}

We next consider the set $$\Omega=\{(\dot\Lambda, \ddot\Lambda, \Lambda)|\dot\Lambda\in P_+^{k}, \ddot\Lambda\in P_+^{l}, \Lambda\in P_+^{k+l}~\text{such that }~\dot{\Lambda}+\ddot{\Lambda}-\Lambda\in Q\}.$$
By Propositions \ref{decomp}, \ref{iden}, we have $\dot\Lambda^{(h)}+ \ddot\Lambda^{(h)}- \Lambda^{(h)}\in Q$ if  $\dot{\Lambda}+\ddot{\Lambda}-\Lambda\in Q$ and $h+Q^{\vee}\in P^{\vee}/Q^{\vee}$. Therefore, we have  an action of $P^{\vee}/Q^{\vee}$ on $\Omega$ defined as follows:
\begin{align*}
\pi: P^{\vee}/Q^{\vee}\times \Omega&\to \Omega\\
(h+Q^{\vee}, (\dot\Lambda, \ddot\Lambda, \Lambda))&\mapsto (\dot\Lambda^{(h)}, \ddot\Lambda^{(h)}, \Lambda^{(h)}).
\end{align*}
For any $(\dot\Lambda, \ddot\Lambda, \Lambda)\in \Omega$, we define the stabilizer $(P^{\vee}/Q^{\vee})_{(\dot\Lambda, \ddot\Lambda, \Lambda)}$ of $(\dot\Lambda, \ddot\Lambda, \Lambda)\in \Omega$ by $$(P^{\vee}/Q^{\vee})_{(\dot\Lambda, \ddot\Lambda, \Lambda)}=\{h+Q^{\vee}\in P^{\vee}/Q^{\vee}|(\dot\Lambda, \ddot\Lambda, \Lambda)=(\dot\Lambda^{(h)}, \ddot\Lambda^{(h)}, \Lambda^{(h)})\}.$$ We use $\Omega/(P^{\vee}/Q^{\vee})$ to denote the set of orbits of $\Omega$ under the action of $P^{\vee}/Q^{\vee}$. For any  $(\dot\Lambda, \ddot\Lambda, \Lambda)\in \Omega$, we denote the orbit of $(\dot\Lambda, \ddot\Lambda, \Lambda)$ by $[\dot\Lambda, \ddot\Lambda, \Lambda]$. Then we have
\begin{lemma}\label{orbit}
Assume that for any  $(\dot\Lambda, \ddot\Lambda, \Lambda)\in \Omega$, the stabilizer $(P^{\vee}/Q^{\vee})_{(\dot\Lambda, \ddot\Lambda, \Lambda)}$ is trivial. Then for any $(\dot\Lambda, \ddot\Lambda, \Lambda)\in \Omega$, the orbit $[\dot\Lambda, \ddot\Lambda, \Lambda]$ of $(\dot\Lambda, \ddot\Lambda, \Lambda)$ has $|P^{\vee}/Q^{\vee}|$ elements.
\end{lemma}

To classify irreducible modules of $C(L_{\g}(k+l,0),L_{\g}(k,0)\otimes L_{\g}(l,0))$, we also need the following fact which was obtained in Corollary 2.7 of \cite{KW}.
\begin{proposition}\label{S3}
Let $\g$ be a finite dimensional simple Lie algebra and $k$ be a positive integer. Then the sum of $S_{L_{\g}(k, 0), L_{\g}(k, \Lambda)}^2$, where $\Lambda$ runs over a congruence class of $P_+^k$ mod $Q$, is equal to $|P/Q|^{-1}$.
\end{proposition}

We are now ready to prove the main result in this section.
\begin{theorem}\label{classi}
Let $k, l$ be positive integers. Suppose that the vertex operator algebra $C(L_{\g}(k+l,0),L_{\g}(k,0)\otimes L_{\g}(l,0))$ is rational, $C_2$-cofinite and that the stabilizer $(P^{\vee}/Q^{\vee})_{(\dot\Lambda, \ddot\Lambda, \Lambda)}$ of $(\dot\Lambda, \ddot\Lambda, \Lambda)\in \Omega$ is trivial for any  $(\dot\Lambda, \ddot\Lambda, \Lambda)\in \Omega$. Then \\
(1) For any $(\dot\Lambda, \ddot\Lambda, \Lambda)\in \Omega$, $M_{\dot{\Lambda}, \ddot{\Lambda}}^{\Lambda}$ is an irreducible $C(L_{\g}(k+l,0),L_{\g}(k,0)\otimes L_{\g}(l,0))$-module.\\
(2) $\{M_{\dot{\Lambda}, \ddot{\Lambda}}^{\Lambda}|[\dot\Lambda, \ddot\Lambda, \Lambda]\in \Omega/(P^{\vee}/Q^{\vee})\}$ is the complete list of inequivalent irreducible modules of  $C(L_{\g}(k+l,0),L_{\g}(k,0)\otimes L_{\g}(l,0))$.
\end{theorem}
\pf It is enough to show that $\{M_{\dot{\Lambda}, \ddot{\Lambda}}^{\Lambda}|[\dot\Lambda, \ddot\Lambda, \Lambda]\in \Omega/(P^{\vee}/Q^{\vee})\}$ satisfies the condition (ii) of Theorem \ref{general}. First, for any $\dot\Lambda\in P_+^{k}, \ddot\Lambda\in P_+^{l}, \Lambda\in P_+^{k+l}$,  it follows from  Proposition \ref{iden} that $M_{\dot{\Lambda}, \ddot{\Lambda}}^{\Lambda}$ is isomorphic to one of modules in $\{M_{\dot{\Lambda}, \ddot{\Lambda}}^{\Lambda}|[\dot\Lambda, \ddot\Lambda, \Lambda]\in \Omega/(P^{\vee}/Q^{\vee})\}$.

Moreover, by Propositions \ref{iden1}, \ref{S3}, Theorems \ref{miden}, \ref{qdimc} and Lemma \ref{orbit}, we have
\begin{align*}
&\sum_{[\dot\Lambda, \ddot\Lambda, \Lambda]\in \Omega/(P^{\vee}/Q^{\vee})}(\qdim_{C(L_{\g}(k+l,0),L_{\g}(k,0)\otimes L_{\g}(l,0))} M_{\dot{\Lambda}, \ddot{\Lambda}}^{\Lambda})^2\\
&=\sum_{(\dot\Lambda, \ddot\Lambda, \Lambda)\in \Omega}\frac{1}{|P^{\vee}/Q^{\vee}|}(\qdim_{C(L_{\g}(k+l,0),L_{\g}(k,0)\otimes L_{\g}(l,0))} M_{\dot{\Lambda}, \ddot{\Lambda}}^{\Lambda})^2\\
&=\sum_{(\dot\Lambda, \ddot\Lambda, \Lambda)\in \Omega}\frac{1}{|P/Q|} \left(\frac{S_{L_{\g}(k, 0), L_{\g}(k, \dot\Lambda)}S_{L_{\g}(l, 0), L_{\g}(l, \ddot\Lambda)}S_{L_{\g}(k+l, 0), L_{\g}(k+l, \Lambda)}}{S_{L_{\g}(k, 0), L_{\g}(k, 0)}S_{L_{\g}(l, 0), L_{\g}(l, 0)}S_{L_{\g}(k+l, 0), L_{\g}(k+l, 0)} }\right)^2\\
&=\sum_{\dot\Lambda\in P^k_+, \ddot\Lambda\in P^l_+}\frac{1}{|P/Q|^2} \left(\frac{S_{L_{\g}(k, 0), L_{\g}(k, \dot\Lambda)}S_{L_{\g}(l, 0), L_{\g}(l, \ddot\Lambda)}}{S_{L_{\g}(k, 0), L_{\g}(k, 0)}S_{L_{\g}(l, 0), L_{\g}(l, 0)}S_{L_{\g}(k+l, 0), L_{\g}(k+l, 0)} }\right)^2\\
&=\frac{1}{|P/Q|^2} \left(\frac{1}{S_{L_{\g}(k, 0), L_{\g}(k, 0)}S_{L_{\g}(l, 0), L_{\g}(l, 0)}S_{L_{\g}(k+l, 0), L_{\g}(k+l, 0)} }\right)^2\\
&=\Glob C(L_{\g}(k+l,0),L_{\g}(k,0)\otimes L_{\g}(l,0)).
\end{align*}
Therefore, by Theorem \ref{general}, $M_{\dot{\Lambda}, \ddot{\Lambda}}^{\Lambda}$  is an irreducible module of $C(L_{\g}(k+l,0),L_{\g}(k,0)\otimes L_{\g}(l,0))$ for any $(\dot\Lambda, \ddot\Lambda, \Lambda)\in \Omega$. Moreover, $\{M_{\dot{\Lambda}, \ddot{\Lambda}}^{\Lambda}|[\dot\Lambda, \ddot\Lambda, \Lambda]\in \Omega/(P^{\vee}/Q^{\vee})\}$ is the complete list of inequivalent irreducible modules of  $C(L_{\g}(k+l,0),L_{\g}(k,0)\otimes L_{\g}(l,0))$.
\qed
\begin{remark}
The multiplicity spaces of $L_{\g}(k+l, \Lambda)$ in $L_{\g}(k, \dot\Lambda)\otimes L_{\g}(l, \ddot\Lambda)$, $\dot\Lambda\in P_+^{k}, \ddot\Lambda\in P_+^{l}, \Lambda\in P_+^{k+l}$, have also been studied in the framework of conformal nets \cite{X}. Many important results have been established in \cite{X}. In particular, our results in Theorem \ref{classi} are motivated by Theorem 4.6 of \cite{X}.
\end{remark}
\section{Trace functions of diagonal coset vertex operator algebras}
\def\theequation{5.\arabic{equation}}
\setcounter{equation}{0}
In this section, our goal is to obtain an explicit modular transformation formula of traces functions of diagonal coset vertex operator algebras.
\subsection{Convergence of one-point theta functions of affine vertex operator algebras} In this subsection, we show that  one-point theta functions of affine vertex operator algebras are convergent, then we can use Theorem \ref{theta2} to obtain the modular transformation formula of traces functions of diagonal coset vertex operator algebras.

 Let $\g$ be a finite dimensional simple Lie algebra and  $k, l$ be positive integers. Then the vertex operator algebra $L_{\g}(k,0)\otimes L_{\g}(l,0)$ is strongly regular. By the discussion in subsection 3.1,  $L_{\g}(k,0)\otimes L_{\g}(l,0)$ has subalgebras $L_{\g}(k+l,0)$ and $C(L_{\g}(k+l,0),L_{\g}(k,0)\otimes L_{\g}(l,0))$ such that $L_{\g}(k+l,0)\otimes C(L_{\g}(k+l,0),L_{\g}(k,0)\otimes L_{\g}(l,0))$ is a vertex subalgebra of $L_{\g}(k,0)\otimes L_{\g}(l,0)$. Moreover, by Proposition \ref{scoset}, $L_{\g}(k+l,0)\otimes C(L_{\g}(k+l,0),L_{\g}(k,0)\otimes L_{\g}(l,0))$ and $L_{\g}(k,0)\otimes L_{\g}(l,0)$ have the same Virasoro vector.  In this section, we always assume that $C(L_{\g}(k+l,0),L_{\g}(k,0)\otimes L_{\g}(l,0))$ is rational and $C_2$-cofinite. Then it follows from Theorem \ref{crational} that $C(L_{\g}(k+l,0),L_{\g}(k,0)\otimes L_{\g}(l,0))$ is strongly regular.

In the following, we use $\w, \w^a$ and $\w^d$ to denote the Virasoro vectors of $L_{\g}(k,0)\otimes L_{\g}(l,0)$, $L_{\g}(k+l,0)$ and $C(L_{\g}(k+l,0),L_{\g}(k,0)\otimes L_{\g}(l,0))$, respectively. We then use $L(n)$, $L^a(n)$ and $L^d(n)$ to denote the component operators of $Y(\w, z)$, $Y(\w^a, z)$ and $Y(\w^d, z)$ respectively for $n\in \Z$. We also denote the central charges of $L_{\g}(k,0)\otimes L_{\g}(l,0)$, $L_{\g}(k+l,0)$ and $C(L_{\g}(k+l,0),L_{\g}(k,0)\otimes L_{\g}(l,0))$ by $c$, $c^a$ and $c^d$, respectively. By Proposition \ref{scoset}, we have $\w=\w^a+\w^d$, $L(0)=L^a(0)+L^d(0)$ and $c=c^a+c^d$.

We now define the one-point theta function of $L_{\g}(k,0)\otimes L_{\g}(l,0))$ as in subsection 2.2. We use $\h^a$ to denote the subset $\{x(-1)\1\otimes \1+\1\otimes x(-1)\1|x\in \h\}$ of $L_{\g}(k,0)\otimes L_{\g}(l,0))$. Note that $\h^a$ is actually a Cartan subalgebra of the weight one Lie algebra $L_{\g}(k+l,0)_1$ of the diagonal subalgebra $L_{\g}(k+l,0)$. For any irreducible $L_{\g}(k,0)\otimes L_{\g}(l,0)$-module $L_{\g}(k,\dot \Lambda)\otimes L_{\g}(l,\ddot\Lambda)$, we  define
$$\chi_{\dot\Lambda, \ddot\Lambda}(w, u, q)=\tr_{L_{\g}(k,\dot \Lambda)\otimes L_{\g}(l,\ddot\Lambda)}o(w)e^{2\pi \sqrt{-1}o(u)}q^{L(0)-c/24},$$
for $w\in C(L_{\g}(k+l,0),L_{\g}(k,0)\otimes L_{\g}(l,0))$, $u\in \h^a$.

To study the convergence of $\chi_{\dot\Lambda, \ddot\Lambda}(w, u, q)$, we need to recall from \cite{DKR} some facts about generalized theta functions of $L_{\g}(k+l,0)$. For any irreducible $L_{\g}(k+l,0)$-module $L_{\g}(k+l,\Lambda)$ and $u\in \h^a$, we define
\begin{align}\label{gtheta}
\chi_{\Lambda}(u, q)=\tr_{L_{\g}(k+l,\Lambda)}e^{2\pi\sqrt{-1}u(0)}q^{L^a(0)-c^a/24}.
\end{align}
Then we have the following results (see Theorems 3.1, 4.2 and Proposition 3.2 of \cite{DKR}).
\begin{theorem}\label{miden1}
(1) For any irreducible $L_{\g}(k+l,0)$-module $L_{\g}(k+l,\Lambda)$ and $u\in \h^a$,
$\chi_{\Lambda}(u, q)$  converges to a holomorphic function $\chi_{\Lambda}(u, \tau)$ in $\mathcal{H}$ with $q=e^{2\pi \sqrt{-1}\tau}$.\\
(2) For any  $u\in\h^a$,
$$\chi_{\Lambda}(\frac{u}{\tau},\frac{-1}{\tau})=e^{\pi\sqrt{-1}(k+l)\langle u, u\rangle/\tau}\sum_{\Lambda'\in P^{k+l}_+}S_{L_{\g}(k+l, \Lambda), L_{\g}(k+l, \Lambda')}\chi_{\Lambda'}(u, \tau),$$ where $S_{L_{\g}(k+l, \Lambda), L_{\g}(k+l, \Lambda')}$ is determined by the formula in  Theorem \ref{miden}.
\end{theorem}

We are now ready to prove the main result in this subsection.
\begin{theorem}\label{conv}
Let $\g$ be a finite dimensional simple Lie algebra and  $k, l$ be positive integers. Suppose that $C(L_{\g}(k+l,0),L_{\g}(k,0)\otimes L_{\g}(l,0))$ is rational and $C_2$-cofinite. Then for any irreducible $L_{\g}(k,0)\otimes L_{\g}(l,0)$-module $L_{\g}(k,\dot \Lambda)\otimes L_{\g}(l,\ddot\Lambda)$,
$\chi_{\dot\Lambda, \ddot\Lambda}(w, u, q)$ converges to a holomorphic function $\chi_{\dot\Lambda, \ddot\Lambda}(w, u, \tau)$ in $\mathcal{H}$
 with $q=e^{2\pi \sqrt{-1}\tau}$ for any $w\in C(L_{\g}(k+l,0),L_{\g}(k,0)\otimes L_{\g}(l,0))$ and $u\in \h^a$.
\end{theorem}
\pf By Proposition \ref{decomp}, $L_{\g}(k, \dot{\Lambda})\otimes L_{\g}(l, \ddot{\Lambda})$ viewed as an $L_{\g}(k+l, 0)\otimes C(L_{\g}(k+l,0),L_{\g}(k,0)\otimes L_{\g}(l,0))$-module has the following decomposition
$$L_{\g}(k, \dot{\Lambda})\otimes L_{\g}(l, \ddot{\Lambda})=\oplus_{\Lambda \in P_+^{k+l}; \dot{\Lambda}+\ddot{\Lambda}-\Lambda\in Q}L_{\g}(k+l, \Lambda)\otimes M_{\dot{\Lambda}, \ddot{\Lambda}}^{\Lambda}.$$
Thus, for any $w\in C(L_{\g}(k+l,0),L_{\g}(k,0)\otimes L_{\g}(l,0))$ and $u\in \h^a$, we have
\begin{align}\label{diden}
\chi_{\dot\Lambda, \ddot\Lambda}(w, u, q)&=\tr_{L_{\g}(k,\dot \Lambda)\otimes L_{\g}(l,\ddot\Lambda)}o(w)e^{2\pi \sqrt{-1}o(u)}q^{L(0)-c/24}\notag\\
&=\sum_{\Lambda \in P_+^{k+l}; \dot{\Lambda}+\ddot{\Lambda}-\Lambda\in Q}\tr_{L_{\g}(k+l, \Lambda)\otimes M_{\dot{\Lambda}, \ddot{\Lambda}}^{\Lambda}}o(w)e^{2\pi \sqrt{-1}o(u)}q^{L^a(0)+L^d(0)-c^a/24-c^d/24}\notag\\
&=\sum_{\Lambda \in P_+^{k+l}; \dot{\Lambda}+\ddot{\Lambda}-\Lambda\in Q}\tr_{ M_{\dot{\Lambda}, \ddot{\Lambda}}^{\Lambda}}o(w)q^{L^d(0)-c^d/24}\tr_{L_{\g}(k+l, \Lambda)}e^{2\pi \sqrt{-1}o(u)}q^{L^a(0)-c^a/24}\notag\\
&=\sum_{\Lambda \in P_+^{k+l}; \dot{\Lambda}+\ddot{\Lambda}-\Lambda\in Q}Z_{ M_{\dot{\Lambda}, \ddot{\Lambda}}^{\Lambda}}(w, q)\chi_{\Lambda}(u, q).
\end{align}
By assumption, $C(L_{\g}(k+l,0),L_{\g}(k,0)\otimes L_{\g}(l,0))$ is rational and $C_2$-cofinite, hence $Z_{ M_{\dot{\Lambda}, \ddot{\Lambda}}^{\Lambda}}(w, q)$ converges to a holomorphic function $Z_{ M_{\dot{\Lambda}, \ddot{\Lambda}}^{\Lambda}}(w, \tau)$ in $\mathcal H$ with $q=e^{2\pi \sqrt{-1}\tau}$. By Theorem \ref{miden1}, $\chi_{\Lambda}(u, q)$  converges to a holomorphic function $\chi_{\Lambda}(u, \tau)$ in $\mathcal{H}$ with $q=e^{2\pi \sqrt{-1}\tau}$. Therefore, $\chi_{\dot\Lambda, \ddot\Lambda}(w, u, q)$  converges to a holomorphic function $\chi_{\dot\Lambda, \ddot\Lambda}(w, u, \tau)$ in $\mathcal{H}$ with $q=e^{2\pi \sqrt{-1}\tau}$.
\qed
\subsection{Modular invariance properties of $C(L_{\g}(k+l,0),L_{\g}(k,0)\otimes L_{\g}(l,0))$} In this subsection, we shall prove the modular transformation formula of traces functions of diagonal coset vertex operator algebras under certain assumption. 
We shall need the following result which was essentially established in Proposition 4.3 of \cite{KP}.
\begin{proposition}\label{inde}
Let $\mathcal O(\mathcal H)$ be the ring of holomorphic functions of $\tau\in \mathcal H$. Then the functions
$\{\chi_{\Lambda}(u, \tau)|\Lambda\in P^{k+l}_+\}$ defined by (\ref{gtheta}) are linearly independent over $\mathcal O(\mathcal H)$.
\end{proposition}

We are now ready to prove the main result in this section.
\begin{theorem}\label{main}
Let $\g$ be a finite dimensional simple Lie algebra and  $k, l$ be positive integers. Suppose  that $C(L_{\g}(k+l,0),L_{\g}(k,0)\otimes L_{\g}(l,0))$ is rational and $C_2$-cofinite. Then for any homogeneous $w\in C(L_{\g}(k+l,0),L_{\g}(k,0)\otimes L_{\g}(l,0))$ and $(\dot\Lambda, \ddot\Lambda, \Lambda)\in \Omega$, we have
\begin{align*}
&Z_{M_{\dot{\Lambda}, \ddot{\Lambda}}^{\Lambda}}(w, \frac{-1}{\tau})=\tau^{\wt [w]}\\
&\cdot\sum_{(\dot\Lambda', \ddot\Lambda', \Lambda')\in \Omega} S_{L_{\g}(k, \dot\Lambda), L_{\g}(k, \dot\Lambda')}S_{L_{\g}(l, \ddot\Lambda), L_{\g}(l, \ddot\Lambda')}\overline{S_{L_{\g}(k+l, \Lambda), L_{\g}(k+l, \Lambda')}}Z_{M_{\dot{\Lambda'}, \ddot{\Lambda'}}^{\Lambda'}}(w, \tau),
 \end{align*}
 where $S_{L_{\g}(k, \dot\Lambda), L_{\g}(k, \dot\Lambda')}$, $S_{L_{\g}(l, \ddot\Lambda), L_{\g}(l, \ddot\Lambda')}$ and $S_{L_{\g}(k+l, \Lambda), L_{\g}(k+l, \Lambda')}$ are determined by the formula in  Theorem \ref{miden}.
\end{theorem}
\pf By Theorem \ref{conv}, $\chi_{\dot\Lambda, \ddot\Lambda}(w, u, q)$ converges to a holomorphic function $\chi_{\dot\Lambda, \ddot\Lambda}(w, u, \tau)$ in $\mathcal{H}$
 with $q=e^{2\pi \sqrt{-1}\tau}$ for any homogeneous $w\in C(L_{\g}(k+l,0),L_{\g}(k,0)\otimes L_{\g}(l,0))$ and $u\in \h^a$. Therefore, by Theorem \ref{theta2}, we have

 \begin{align}\label{tiden}
 \chi_{\dot\Lambda, \ddot\Lambda}(w, \frac{u}{\tau}, \frac{-1}{\tau})&=\tau^{\wt [w]}e^{\pi\sqrt{-1}(k+l)(\langle u, u\rangle/\tau)}\notag\\
 &\sum_{\dot\Lambda'\in P^k_+, \ddot\Lambda'\in P^l_+} S_{L_{\g}(k, \dot\Lambda), L_{\g}(k, \dot\Lambda')}S_{L_{\g}(l, \ddot\Lambda), L_{\g}(l, \ddot\Lambda')}\chi_{\dot\Lambda', \ddot\Lambda'}(w, u,\tau),
  \end{align}
 where $S_{L_{\g}(k, \dot\Lambda), L_{\g}(k, \dot\Lambda')}$ and $S_{L_{\g}(l, \ddot\Lambda), L_{\g}(l, \ddot\Lambda')}$ are determined by the formula in  Theorem \ref{miden}.

 We next consider column vectors $$\overrightarrow{\chi(w, u, \tau)}=(\chi_{\dot\Lambda, \ddot\Lambda}(w, u, \tau))_{\dot\Lambda\in P^{k}_+, \ddot\Lambda\in P^l_+},~~~\overrightarrow{\chi(u, \tau)}=(\chi_{\Lambda}(u, \tau))_{\Lambda\in P^{k+l}_+}.$$
 Then, by the formula (\ref{tiden}), we have
 $$\overrightarrow{\chi(w, \frac{u}{\tau}, \frac{-1}{\tau})}=\tau^{\wt [w]}e^{\pi\sqrt{-1}(k+l)(\langle u, u\rangle/\tau)}S\cdot \overrightarrow{\chi(w, u, \tau)},$$ where $S$ is the matrix $(S_{\dot\Lambda, \ddot\Lambda; \dot\Lambda', \ddot\Lambda'})_{\dot\Lambda\in P^k_+, \ddot\Lambda\in P^l_+; \dot\Lambda'\in P^k_+, \ddot\Lambda'\in P^l_+}$ such that $$S_{\dot\Lambda, \ddot\Lambda; \dot\Lambda', \ddot\Lambda'}=S_{L_{\g}(k, \dot\Lambda), L_{\g}(k, \dot\Lambda')}S_{L_{\g}(l, \ddot\Lambda), L_{\g}(l, \ddot\Lambda')}.$$
 On the other hand, by Theorem \ref{miden1}, we have
 $$\overrightarrow{\chi(\frac{u}{\tau}, \frac{-1}{\tau})}=e^{\pi\sqrt{-1}(k+l)\langle u, u\rangle/\tau}S^a\cdot\overrightarrow{\chi(u, \tau)},$$
 where $S^a$ is the matrix $(S^a_{\Lambda, \Lambda'})_{\Lambda\in P^{k+l}_+, \Lambda'\in P^{k+l}_+}$ such that $S^a_{\Lambda, \Lambda'}=S_{L_{\g}(k+l, \Lambda), L_{\g}(k+l, \Lambda')}$.

 We also consider the matrix $$Z(w, \tau)=(Z(w, \tau)_{\dot\Lambda,\ddot\Lambda;\Lambda})_{\dot\Lambda\in P^{k}_+, \ddot\Lambda\in P^l_+; \Lambda\in P^{k+l}_+}$$ such that $Z(w, \tau)_{\dot\Lambda,\ddot\Lambda;\Lambda}=Z_{M_{\dot{\Lambda}, \ddot{\Lambda}}^{\Lambda}}(w, \tau)$ if $(\dot\Lambda, \ddot\Lambda, \Lambda)\in \Omega$ and $Z(w, \tau)_{\dot\Lambda,\ddot\Lambda;\Lambda}=0$ otherwise. Therefore, by the formula (\ref{diden}), we have
 $$\overrightarrow{\chi(w, u, \tau)}=Z(w, \tau)\overrightarrow{\chi(u, \tau)}.$$
 Performing the transformation of both sides of matrix $\left(\begin{array}{cc}0 & -1\\ 1 & 0\end{array}\right)$, we see that
 $$\tau^{\wt [w]}e^{\pi\sqrt{-1}(k+l)(\langle u, u\rangle/\tau)}S\cdot \overrightarrow{\chi(w, u, \tau)}=Z(w, \frac{-1}{\tau})e^{\pi\sqrt{-1}(k+l)\langle u, u\rangle/\tau}S^a\cdot\overrightarrow{\chi(u, \tau)}.$$Hence, we have
 $$\tau^{\wt [w]}S\cdot Z(w, \tau)\overrightarrow{\chi(u, \tau)}=Z(w, \frac{-1}{\tau})S^a\cdot\overrightarrow{\chi(u, \tau)}.$$
 By Proposition \ref{inde}, the functions $\{\chi_{\Lambda}(u, \tau)|\Lambda\in P^{k+l}_+\}$ are linearly independent over $\mathcal O(\mathcal H)$. This implies that
 $$\tau^{\wt [w]}S\cdot Z(w, \tau)=Z(w, \frac{-1}{\tau})S^a.$$
 It is known in (13.8.5), (13.8.6) and (13.8.7) of \cite{K} that $S^a$ is symmetric,  unitary  and $(S^{a})^{-1}=\overline{S^a}$. Hence, we have
 $$Z(w, \frac{-1}{\tau})=\tau^{\wt [w]}S\cdot Z(w, \tau)\overline{S^a}.$$
 Comparing the $(\dot\Lambda,\ddot\Lambda;\Lambda)$-entries of both sides, we get
\begin{align*}
&Z_{M_{\dot{\Lambda}, \ddot{\Lambda}}^{\Lambda}}(w, \frac{-1}{\tau})=\tau^{\wt [w]}\\
&\cdot\sum_{(\dot\Lambda', \ddot\Lambda', \Lambda')\in \Omega} S_{L_{\g}(k, \dot\Lambda), L_{\g}(k, \dot\Lambda')}S_{L_{\g}(l, \ddot\Lambda), L_{\g}(l, \ddot\Lambda')}\overline{S_{L_{\g}(k+l, \Lambda), L_{\g}(k+l, \Lambda')}}Z_{M_{\dot{\Lambda'}, \ddot{\Lambda'}}^{\Lambda'}}(w, \tau).
 \end{align*}
This completes the proof.
\qed
\begin{remark}
 In the case that $w$ is the vacuum vector of $C(L_{\g}(k+l,0),L_{\g}(k,0)\otimes L_{\g}(l,0))$, this result has been established in \cite{KW}.
\end{remark}

 We now assume  that $C(L_{\g}(k+l,0),L_{\g}(k,0)\otimes L_{\g}(l,0))$ is rational, $C_2$-cofinite and that the stabilizer $(P^{\vee}/Q^{\vee})_{(\dot\Lambda, \ddot\Lambda, \Lambda)}$ of $(\dot\Lambda, \ddot\Lambda, \Lambda)\in \Omega$ is trivial for any  $(\dot\Lambda, \ddot\Lambda, \Lambda)\in \Omega$. By Theorem \ref{classi}, $$\{M_{\dot{\Lambda}, \ddot{\Lambda}}^{\Lambda}|[\dot\Lambda, \ddot\Lambda, \Lambda]\in \Omega/(P^{\vee}/Q^{\vee})\}$$ is the complete list of inequivalent irreducible modules of  $C(L_{\g}(k+l,0),L_{\g}(k,0)\otimes L_{\g}(l,0))$.
By Theorem \ref{main}, we have
\begin{corollary}\label{smatrix}
Let $\g$ be a finite dimensional simple Lie algebra and  $k, l$ be positive integers. Suppose  that $C(L_{\g}(k+l,0),L_{\g}(k,0)\otimes L_{\g}(l,0))$ is rational, $C_2$-cofinite and that the stabilizer $(P^{\vee}/Q^{\vee})_{(\dot\Lambda, \ddot\Lambda, \Lambda)}$ of $(\dot\Lambda, \ddot\Lambda, \Lambda)\in \Omega$ is trivial for any  $(\dot\Lambda, \ddot\Lambda, \Lambda)\in \Omega$. Then
for any $(\dot\Lambda, \ddot\Lambda, \Lambda), (\dot\Lambda', \ddot\Lambda', \Lambda')\in \Omega$, we have
\begin{align*}
S_{ M_{\dot{\Lambda}, \ddot{\Lambda}}^{\Lambda},  M_{\dot{\Lambda'}, \ddot{\Lambda'}}^{\Lambda'}}=|P/Q|S_{L_{\g}(k, \dot\Lambda), L_{\g}(k, \dot\Lambda')}S_{L_{\g}(l, \ddot\Lambda), L_{\g}(l, \ddot\Lambda')}\overline{S_{L_{\g}(k+l, \Lambda), L_{\g}(k+l, \Lambda')}}.
 \end{align*}
\end{corollary}
\section{Fusion rules of $C(L_{E_8}(k+2,0), L_{E_8}(k,0)\otimes L_{E_8}(2,0))$}
\def\theequation{6.\arabic{equation}}
\setcounter{equation}{0}
As an application of the results obtained in Section 5, we determine the fusion rules of the coset vertex operator algebra $C(L_{E_8}(k+2,0), L_{E_8}(k,0)\otimes L_{E_8}(2,0))$ in this section. First, we recall some facts about fusion
rules from \cite{FHL}.
Let $V$ be a vertex operator algebra and $M^1$, $M^2$, $M^3$ be admissible $V$-modules. An {\em intertwining
operator} $\mathcal {Y}$ of type $\left(\begin{tabular}{c}
$M^3$\\
$M^1$ $M^2$\\
\end{tabular}\right)$ is a linear map
\begin{align*}
\mathcal
{Y}: M^1&\rightarrow \Hom(M^2, M^3)\{z\},\\
 w^1&\mapsto\mathcal {Y}(w^1, z) = \sum_{n\in \mathbb{C}}{w_n^1z^{-n-1}}
\end{align*}
satisfying a number of conditions (cf. \cite{FHL}).  We use $\mathcal{I}_{M^1,M^2}^{M^3}$ to denote the vector space of intertwining operators of type $\left(\begin{tabular}{c}
$M^3$\\
$M^1$ $M^2$\\
\end{tabular}\right)$.  If $V$ is a strongly regular vertex operator algebra and $M^0,...,M^p$  are all the inequivalent irreducible $V$-modules,  we define the \emph{fusion rules} to be the formal product rules
 \begin{align*}
 M^i\times M^j=\sum_{0\leq k\leq p} N_{M^i,M^j}^{M^k}M^k,
 \end{align*}
 where $N_{M^i,M^j}^{M^k}$ denotes the dimension  of $\mathcal{I}_{M^i,M^j}^{M^k}$. To determine the fusion rules, we need the Verlinde formula proved in  \cite{H}.
 \begin{theorem}\label{verlinde}
 Let $V$ be a strongly regular vertex operator algebra,  $M^0=V, M^1,...,M^p$  be all the inequivalent irreducible $V$-modules.  Then
 $$N_{M^i,M^j}^{M^k}=\sum_{a=0}^{p}\frac{S_{j, a}S_{i, a}(S^{-1})_{k, a}}{S_{0,a}}.$$
 \end{theorem}

\vskip.25cm
We next recall some facts about  $C(L_{E_8}(k+2, 0), L_{E_8}(k, 0)\otimes L_{E_8}(2, 0))$. First, we have the following result proved in Theorem 3.10 of \cite{Lin}.
\begin{theorem}\label{main1}
Let $k$ be a positive integer. Then the coset vertex operator algebra $C(L_{E_8}(k+2, 0), L_{E_8}(k, 0)\otimes L_{E_8}(2, 0))$ is rational and $C_2$-cofinite.
\end{theorem}

Note that if $\g=E_8$ then $|P/Q|=1$. Then the following result about classification of irreducible modules of  $C(L_{E_8}(k+2,0), L_{E_8}(k,0)\otimes L_{E_8}(2,0))$ has been obtained in Theorem 5.3 of \cite{Lin}.
 \begin{theorem}\label{main2}
 Let $k$ be a positive integer. Then  $$\{ M_{\dot{\Lambda}, \ddot{\Lambda}}^{\Lambda}|\dot\Lambda\in P_+^{k}, \ddot\Lambda\in P_+^{2}, \Lambda\in P_+^{k+2}\}$$
 are all the inequivalent irreducible $C(L_{E_8}(k+2,0), L_{E_8}(k,0)\otimes L_{E_8}(2,0))$-modules.
 \end{theorem}
 We are now ready to determine the fusion rules of $C(L_{E_8}(k+2,0), L_{E_8}(k,0)\otimes L_{E_8}(2,0))$.
 \begin{theorem}\label{fusion}
 Let $k$ be a positive integer. Then we have for any $\dot\lambda, \dot\mu, \dot\nu\in P^{k}_+$, $\ddot\lambda, \ddot\mu, \ddot\nu\in P^{2}_+$ and $\lambda, \mu, \nu\in P^{k+2}_+$,
 $$N_{M_{\dot{\lambda}, \ddot{\lambda}}^{\lambda}, M^{\mu}_{\dot{\mu}, \ddot{\mu}}}^{M_{\dot{\nu}, \ddot{\nu}}^{\nu}}=N_{L_{E_8}(k, \dot{\lambda}), L_{E_8}(k, \dot{\mu})}^{L_{E_8}(k, \dot{\nu})}N_{L_{E_8}(2, \ddot{\lambda}), L_{E_8}(2, \ddot{\mu})}^{L_{E_8}(2, \ddot{\nu})}N_{L_{E_8}(k+2, \lambda), L_{E_8}(k+2, \mu)}^{L_{E_8}(k+2, \nu)}.$$
 \end{theorem}
 \pf First, by Theorems \ref{crational}, \ref{main1}, we know that $C(L_{E_8}(k+2,0), L_{E_8}(k,0)\otimes L_{E_8}(2,0))$ is a strongly regular vertex operator algebra. In the following, for any $\dot\Lambda, \dot\Lambda'\in P^{k}_+$, $\ddot\Lambda, \ddot\Lambda'\in P^{2}_+$ and $\Lambda, \Lambda'\in P^{k+2}_+$, we use $S_{\dot\Lambda, \dot\Lambda'}$, $S_{ \ddot\Lambda, \ddot\Lambda'}$ and $S_{\Lambda, \Lambda'}$ to denote $S_{L_{\g}(k, \dot\Lambda), L_{\g}(k, \dot\Lambda')}$, $S_{L_{\g}(2, \ddot\Lambda), L_{\g}(2, \ddot\Lambda')}$ and $S_{L_{\g}(k+2, \Lambda), L_{\g}(k+2, \Lambda')}$, respectively. Then, by Corollary \ref{smatrix}, for any $(\dot\Lambda, \ddot\Lambda, \Lambda), (\dot\Lambda', \ddot\Lambda', \Lambda')\in \Omega$, we have
\begin{align*}
S_{ M_{\dot{\Lambda}, \ddot{\Lambda}}^{\Lambda},  M_{\dot{\Lambda'}, \ddot{\Lambda'}}^{\Lambda'}}=S_{\dot\Lambda,  \dot\Lambda'}S_{ \ddot\Lambda, \ddot\Lambda'}\overline{ S_{\Lambda, \Lambda'}}.
 \end{align*}
 Therefore, by Theorems \ref{verlinde}, \ref{main2}, and Corollary 3.6 of \cite{DLN}, we have for any $\dot\lambda, \dot\mu, \dot\nu\in P^{k}_+$, $\ddot\lambda, \ddot\mu, \ddot\nu\in P^{2}_+$ and $\lambda, \mu, \nu\in P^{k+2}_+$,
\begin{align*}
 &N_{M_{\dot{\lambda}, \ddot{\lambda}}^{\lambda}, M^{\mu}_{\dot{\mu}, \ddot{\mu}}}^{M_{\dot{\nu}, \ddot{\nu}}^{\nu}}\\
 &=\sum_{\dot\Lambda\in P_+^{k}, \ddot\Lambda\in P_+^{2}, \Lambda\in P_+^{k+2}}\frac{S_{M_{\dot{\lambda},\ddot{\lambda}}^{\lambda}, M_{\dot{\Lambda},\ddot{\Lambda}}^{\Lambda}}
 S_{M^{\mu}_{\dot{\mu},\ddot{\mu}},M_{\dot{\Lambda}, \ddot{\Lambda}}^{\Lambda}}\overline{S_{M_{\dot{\nu}, \ddot{\nu}}^{\nu}, M_{\dot{\Lambda}, \ddot{\Lambda}}^{\Lambda}}}}{S_{C(L_{\g}(k+2,0),L_{\g}(k,0)\otimes L_{\g}(2,0)), M_{\dot{\Lambda}, \ddot{\Lambda}}^{\Lambda}}}\\
 &=\sum_{\dot\Lambda\in P_+^{k}, \ddot\Lambda\in P_+^{2}, \Lambda\in P_+^{k+2}}\frac{S_{\dot\lambda,  \dot\Lambda}S_{ \ddot\lambda, \ddot\Lambda}\overline{ S_{\lambda, \Lambda}}S_{\dot\mu,  \dot\Lambda}S_{ \ddot\mu, \ddot\Lambda}\overline{ S_{\mu, \Lambda}}\overline{S_{\dot\nu,  \dot\Lambda}S_{ \ddot\nu, \ddot\Lambda}\overline{ S_{\nu, \Lambda}}}}{S_{0,  \dot\Lambda}S_{ 0, \ddot\Lambda}\overline{ S_{0, \Lambda}}}\\
 &=\sum_{\dot\Lambda\in P_+^{k}}\frac{S_{\dot\lambda,  \dot\Lambda}S_{\dot\mu,  \dot\Lambda}\overline{S_{\dot\nu,  \dot\Lambda}}}{S_{0,  \dot\Lambda}}\sum_{\ddot\Lambda\in P_+^{2}}\frac{S_{\ddot\lambda,  \ddot\Lambda}S_{\ddot\mu,  \ddot\Lambda}\overline{S_{\ddot\nu,  \ddot\Lambda}}}{S_{0,  \ddot\Lambda}}\sum_{\Lambda\in P_+^{k+2}}\frac{S_{\lambda,  \Lambda}S_{\mu,  \Lambda}\overline{S_{\nu,  \Lambda}}}{S_{0,  \Lambda}}\\
 &=N_{L_{E_8}(k, \dot{\lambda}), L_{E_8}(k, \dot{\mu})}^{L_{E_8}(k, \dot{\nu})}N_{L_{E_8}(2, \ddot{\lambda}), L_{E_8}(2, \ddot{\mu})}^{L_{E_8}(2, \ddot{\nu})}N_{L_{E_8}(k+2, \lambda), L_{E_8}(k+2, \mu)}^{L_{E_8}(k+2, \nu)}.
 \end{align*}
 This completes the proof.
 \qed
 \begin{remark}
 The similar formula of the fusion rules of the vertex operator algebra $C(L_{sl_2}(k+2,0), L_{sl_2}(k,0)\otimes L_{sl_2}(2,0))$ has been claimed in Remark 6.4 of \cite{CFL}.
 \end{remark}

{\bf Acknowledgement}. The author wish to thank Matthew Krauel for explaining the results in \cite{Kra}.


\begin{thebibliography}{ABCD}
\bibitem{ACL}
Arakawa, T., Creutzig,  T.,    Linshaw, A.: W-algebras as coset vertex algebras. {\em Invent. Math.} {\bf 218}, 145-195 (2019)

\bibitem{CFL}
Creutzig, T., Feigin, B., Linshaw, A.: $N = 4$ superconformal algebras
and diagonal cosets. arXiv:1910.01228.
\bibitem
{DJX}
Dong, C., Jiao, X., Xu, F.: Quantum dimensions and quantum Galois theory. {\em Trans. Amer. Math. Soc.} {\bf 365}, 6441-6469  (2013)
\bibitem
{DKR}
 Dong, C., Kac, V.,  Ren, L.:  Trace functions of the parafermion vertex operator algebras. {\em Adv. Math.} {\bf 348}, 1-17  (2019)

\bibitem
{DLM1} Dong, C., Li, H., Mason, G.:
Twisted representations of vertex operator algebras. {\em Math. Ann.}
{\bf  310}, 571--600  (1998)
\bibitem
{DLM2} Dong, C., Li, H., Mason, G.: Modular invariance
of trace functions in orbifold theory and generalized moonshine. {\em
Comm. Math. Phys.} {\bf 214}, 1-56  (2000)

\bibitem{DLN}
Dong, C., Lin, X., Ng, S.: Congruence property in conformal field theory. {\em Algebra Number Theory} {\bf 9}, 2121-2166  (2015)

\bibitem{DLMa}
Dong, C., Liu, K., Ma, X.: Elliptic genus and vertex operator algebras. {\em Pure Appl. Math. Q.} {\bf 1},  791-815  (2005)

\bibitem
{DM1} Dong, C.,  Mason, G.: Rational vertex operator
algebras and the effective central charge. {\em Internat. Math. Res.
Notices} {\bf 56}, 2989-3008  (2004)

\bibitem
{DM2} Dong, C.,  Mason, G.:  Integrability of $C_2$-cofinite vertex operator algebra, {\em Internat. Math. Res.
Notices} {\bf 2006}, Article ID 80468, 15 pages  (2006)

\bibitem
{FHL} Frenkel, I., Huang, Y.,  Lepowsky, J.: On axiomatic
approaches to vertx operator algebras and modules. {\em Mem. Amer. Math. Soc.} {\bf  104} (1993)
\bibitem
{FZ} Frenkel, I.,  Zhu, Y.: Vertex operator algebras
associated to representations of affine and Virasoro algebra. {\em Duke.
Math. J.} {\bf 66}, 123-168 (1992)
 \bibitem{H}
  Huang, Y.: Vertex operator algebras and the Verlinde conjecture. {\em Commun. Contemp. Math.} {\bf 10}, 103-154 (2008)
\bibitem{JL1}
Jiang, C., Lin, Z.: The commutant of $L_{\hat{sl}_2}(n, 0)$
in the vertex operator algebra $L_{\hat{sl}_2}(1, 0)^{\otimes n}$. {\em Adv. Math.} {\bf 301}, 227-257  (2016)

\bibitem{JL2}
Jiang, C., Lin, Z.: Tensor decomposition, parafermions, level-rank
duality, and reciprocity laws for vertex operator algebras. arXiv:1406.4191.
\bibitem
{K} Kac, V.: Infinite dimensional Lie algebras. Third edition. {\em Cambridge University Press, Cambridge} (1990)
\bibitem
{KP}
Kac, V., Peterson, D.: Infinite-dimensional Lie algebras, theta functions and modular forms. {\em Adv. in Math.} {\bf 53},  125-264  (1984)
\bibitem
{KW} Kac, V., Wakimoto, M.: Modular and conformal invariance constraints in representation theory of affine algebras. {\em Adv. Math.} {\bf 70}, 156-236  (1988)
 \bibitem{Kra}
  Krauel, M.: One-point theta functions for vertex operator algebras. {\em J. Algebra} {\bf 481}, 250-272  (2017)

\bibitem
{LL} Lepowsky, J.,  Li, H.: Introduction to vertex operator
algebras and their representations. Progress in Mathematics, 227. Birkh$\ddot{a}$user Boston, Inc., Boston, MA (2004)
\bibitem
{Li}
Li, H.: Symmetric invariant bilinear forms on vertex operator algebras. {\em J. Pure Appl. Algebra} {\bf 96}, 279-297  (1994)

\bibitem
{L3}
 Li, H.: Certain extensions of vertex operator algebras of affine type. {\em Comm. Math. Phys.} {\bf 217},  653-696  (2001)
 \bibitem
{L4}  Li, H.: Local systems of twisted vertex operators, vertex operator superalgebras and twisted modules, in: Moonshine, the Monster, and related topics. {\em Contemp. Math.} {\bf 193}, 203--236. American Mathematical Sociaty, Providence, RI,
   (1995)
 \bibitem{Lin}
 Lin, X.: Quantum dimensions and irreducible modules of some diagonal coset vertex operator algebras. {\em Lett. Math. Phys.} {\bf 110}, 1363-1380  (2020)
\bibitem
{X}
 Xu, F.:  Algebraic coset conformal field theories. {\em Comm. Math. Phys.} {\bf 211},  1-43  (2000)
\bibitem
{Z} Zhu, Y.: Modular invariance of characters of vertex
operator algebras. {\em J. Amer. Math. Soc.} {\bf 9}, 237-302  (1996)
\end{thebibliography}
\end{document}